\documentclass[12pt]{article}

\usepackage{color}
\title{Good reduction of the Brauer--Manin obstruction}
\author{Jean-Louis Colliot-Th\'el\`ene and Alexei N. Skorobogatov}
\date{September 1st, 2010}

\def\g{{\mathfrak g}}

\def\sS{{\mathcal S}}
\def\sZ{{\mathcal Z}}
\def\beq{\begin{equation} \label}

\usepackage{amssymb}
\usepackage{amsmath}
\usepackage{amscd}
\usepackage[arrow]{xy}

\def\sO{{\mathcal O}}

\def\sS{{\mathcal S}}
\def\sX{{\mathcal X}}

\def\qed{{\hfill$\square$}}

\newtheorem{theo}{Theorem}[section]
\newtheorem{prop}[theo]{Proposition}
\newtheorem{lem}[theo]{Lemma}
\newtheorem{cor}[theo]{Corollary}
\newtheorem{defi}[theo]{Definition}

\newcommand{\bthe}{\begin{theo}}
\newcommand{\ble}{\begin{lem}}
\newcommand{\bpr}{\begin{prop}}
\newcommand{\bco}{\begin{cor}}
\newcommand{\bde}{\begin{defi}}
\newcommand{\ethe}{\end{theo}}
\newcommand{\ele}{\end{lem}}
\newcommand{\epr}{\end{prop}}
\newcommand{\eco}{\end{cor}}
\newcommand{\ede}{\end{defi}}

\def\lra{\longrightarrow}
\def\res{{\rm res}}
\def\H{{\rm H}}
\def\NS{{\rm NS}}
\def \Gal {{\rm{Gal}}}
\def \Spec {{\rm{Spec}}}
\def\et{\rm{\acute et}}
\def \Pic {{\rm {Pic}}}
\def \Div {{\rm {Div}}}
\def\ov{\overline}
\def \Br {{\rm{Br}}}
\def \to {{\rightarrow}}
\def \Ker {{\rm Ker}}
\def \Hom {{\rm Hom}}
\def \A {{\mathbb A}}
\def \Q {{\mathbb Q}}
\def \Z {{\mathbb Z}}
\def \F {{\mathbb F}}

\def \R {{\mathbb R}}
\def \X {{\mathcal X}}

\def \G {{\mathbb G}}
\def \P {{\mathbb P}}

\begin{document}

\maketitle

\begin{abstract}
\noindent For a smooth and projective variety over a number field with 
torsion free geometric Picard group and finite transcendental 
Brauer group we show that only the archimedean places, 
the primes of bad reduction and
the primes dividing the order of the transcendental Brauer group 
can turn up in the description of the Brauer--Manin set.
\end{abstract}

\noindent\centerline{AMS Subject Classification: 14F22, 14G05, 11G35, 11G25}

\section*{Introduction}

We assume that the reader is familiar with the basic theory of
the Brauer group \cite{G68} and 
of the Brauer--Manin obstruction, see \cite[\S 3.1]{CS} or \cite[\S 5.2]{S}.

Let us briefly recall some notation.
To any scheme $X$ one associates its (Grothendieck) Brauer group
 $\Br(X)=\H^2_{\et}(X,{\mathbb G}_{m})$. 
Let $k$ be a field,   $\ov k$ be an algebraic closure of $k$, 
$\Gamma={\rm Gal}(\ov k/k)$.
If $X$ is a variety over $k$ we write $\ov X=X\times_k\ov k$.
Recall the standard notation 
$$\Br_0(X)={\rm Im}[\Br(k)\to \Br(X)], \quad
\Br_1(X)={\rm Ker}[\Br(X)\to \Br(\ov X)].$$
The group $\Br_1(X)$ is called the {\it algebraic} Brauer group of $X$,
and the quotient $\Br(X)/\Br_1(X)$ is sometimes called the 
{\it transcendental} Brauer group of $X$, a terminology we
adopt in this note.

Let $k$ be a number field, let $k_v$ be the completion of $k$
at a place $v$, and let $\A_k$ be the ring of ad\`eles of $k$, i.e.
the restricted product of $k_v$ for all places $v$ of $k$, taken 
with respect to the rings of integers $\sO_v\subset k_v$.
Let $\sO$ be the ring of integers of $k$. Given a finite set  $S$
of places of $k$, we let  $\sO_S$
be the subring of $k$ consisting of the elements
that are integral at the primes not in $S$. 

Given a variety $X$ over $k$, the set $ X(\A_k)^\Br \subset X(\A_k)$ 
is the set of ad\`eles of $X$
which are orthogonal to the Brauer group of $X$ with respect to
the Brauer--Manin pairing.

Our aim is to explore what can be proved in 
the direction of the following question.

\bigskip

\noindent{\bf Question 1}. Let $k$ be a number field and
let $S$ be a finite set of places containing the archimedean places.
Let $\X$ be a smooth projective $\sO_S$-scheme 
with geometrically integral fibres.
Let $X/k$ be its generic fibre. Assume that
$\Pic(\ov X)$ is a finitely generated torsion-free abelian group.
Does there exist an open and closed set $Z \subset \prod_{v\in S}X(k_v)$
such that
$$X(\A_k)^\Br=Z\times \prod_{v\not\in S}X(k_v) \hskip2mm ?$$

The assumption on the Picard group is equivalent to the conjunction of
two geometric hypotheses: the coherent cohomology group $\H^1(X,O_{X})$
is zero,
and the geometric N\'eron-Severi group $\NS(\ov X)$ is torsion-free. 

Peter Swinnerton-Dyer asked us this question in connection with 
his recent work \cite{sw}. 
Theorem \ref{gen} below gives an affirmative
answer under the additional hypothesis that
the transcendental Brauer group of $X$ (as defined above)  is finite, 
and $S$ contains all the primes
dividing its order. A criterion that uses 
only a geometric hypothesis on $\ov X$ is given in Corollary \ref{t1}.
For example, the answer is positive whenever $X$ is a surface
of geometric genus $0$. 
In Corollary \ref{t2}, using the results of
\cite{ISZ}, we give a positive answer
to Question~1 for arbitrary smooth diagonal quartic surfaces over $\Q$,
a case which is not covered by Corollary \ref{t1}.
For general varieties, 
 Proposition \ref{f}  relates the finiteness
property of the transcendental Brauer group to the Tate 
conjecture for divisors.

\section{Preliminaries}

Given an abelian group $A$, a positive integer $n$ and a prime number $\ell$,
we denote by $A[n] \subset A$ the group of elements annihilated 
by $n$, and by $A\{\ell\} \subset A$ the group of elements 
annihilated by some power of $\ell$.

The following lemma is well known. 

\ble \label{l1}
Let $X$ be a smooth, projective and geometrically integral 
variety over a number field $k$ such that
$\Pic(\ov X)$ is a finitely generated torsion-free
abelian group, and $\Br(X)/\Br_1(X)$ is finite.
Then $X(\A_k)^\Br$ is open and closed in $X(\A_k)$.
\ele
{\it Proof}  From the spectral sequence 
$$E_{2}^{pq}=\H^p(k,\H^q_{\et}(\ov X,\G_m))\Rightarrow \H^{p+q}_{\et}(X,\G_m)$$
we see that the quotient $\Br_1(X)/\Br_0(X)$ is a subgroup
of $\H^1(k,\Pic(\ov X))$, and this last group is finite
since $\Pic(\ov X)$ is finitely generated and torsion-free. 
Hence $\Br(X)/\Br_0(X)$ is finite.
The sum of local invariants of a given element of $\Br(X)$ is a continuous
function on $X(\A_k)$ with finitely many values, 
and this function is identically zero
if the element is in $\Br_0(X)$. The lemma follows. \qed

\medskip

By the main result of \cite{SZ} the conditions of this lemma 
are satisfied for K3 surfaces.
See Section \ref{Tate} for a discussion of the finiteness
of $\Br(X)/\Br_1(X)$.

\medskip

Let us write $X_v=X\times_k k_v$ 
and, for $v \notin S$, $\sX_v=\sX\times_{\sO_{S}}\sO_v$.

\ble \label{l2}
We keep the assumptions of Lemma $\ref{l1}$. If for every
$v\notin S$ the image of $\Br(X)\to\Br(X_v)$ is contained in the
subgroup generated by the images of $\Br(k_v)$ and $\Br(\sX_v)$,
then the answer to Question $1$ is positive.
\ele
{\it Proof} We may  assume $X(\A_k)^\Br\not=\emptyset$. 
Fix $M=\{M_v\} \in X(\A_k)^\Br$.
From the inclusion 
$$\Br(k) \hookrightarrow  \bigoplus_{{\rm all} \hskip1mm v} \Br(k_{v})$$
provided by class field theory,
and the fact that $X(\A_k)$ is non-empty we conclude that the
natural map
$\Br(k)\rightarrow\Br(X)$ is injective, so that 
$\Br(k)=\Br_0(X)$.
Moreover, the exact sequence of abelian groups
$$0\to\Br(k)\to\Br(X)\to\Br(X)/\Br(k)\to 0$$
has a splitting defined by $M$. Indeed, the injective map
$\Br(k)\to\Br(X)$ has a retraction $r$ sending $A\in \Br(X)$
to the unique element $r(A)\in\Br(k)$ whose image in
$\Br(k_v)$ is $A(M_v)$ for every $v$.
Let us denote by $B$ the kernel of $r:\Br(X)\to\Br(k)$.
An element $A \in \Br(X)$ lies in $B$ if and only if
$A(M_v)=0$ for all $v$.
Then $\Br(X)=\Br(k)\oplus B$, and $B$ projects isomorphically
onto $\Br(X)/\Br(k)$. 

To complete the proof it is enough to show that for any $v\notin S$,
any $A\in B$ and any $N\in X(k_v)$ we have $A(N)=0\in \Br(k_v)$.
By assumption the image of $A$ in $\Br(X_v)$
can be written as a sum $A_{1}+A_{2}$ where $A_{1} \in \Br (\sX_v)$
and $A_{2} \in \Br (k_v)$. 
Since $\sX_v$ is proper over $\Spec(\sO_v)$ we have 
$X(k_v)=\sX_v(\sO_v)$, hence $A_1(N)\in\Br(\sO_v)=0$
for any $N\in X(k_v)$.
We have $A(M_v)=0$ because $A\in B$. It follows that 
the `constant' algebra $A_2$ has value $0$ at $M_v$, hence $A_2=0$.
We have proved that $A(N)=0 \in \Br(k_{v})$, and the lemma follows. \qed

\medskip

The following well known lemma is due to Grothendieck, 
it is given here for the ease of reference.

\ble \label{new1}
Let $X$ be a smooth, projective and geometrically integral variety
over a field of characteristic $0$. Assume that
$\H^1(X,O_{X})=\H^2(X,O_{X})=0$ and that the N\'eron--Severi group
$\NS(\ov X)$ is torsion-free. 

{\rm (i)} Then the groups $\Br(\ov X)$ and
 $\Br(X)/\Br_0(X)$ are finite. 
 
{\rm (ii)} We have $\Br(\ov X)=0$ 
if and only if $\H^3_{\et}(\ov X,\Z_\ell(1))_{\rm tors}=0$
for every prime $\ell$. In this case $\Br_1(X)=\Br(X)$.

{\rm (iii)} If $\dim X=2$, then $\Br(\ov X)=0$
and $\Br_1(X)=\Br(X)$.
\ele
{\it Proof}
By \cite[III.8]{G68} we have a short exact sequence
\begin{equation}
0\lra (\Q/\Z)^{b_{2}-\rho}\lra \Br(\ov X) \lra
\bigoplus_\ell  \H^3_{\et}(\ov X,\Z_\ell(1))_{\rm tors}\lra 0,
\label{br}
\end{equation}
where $b_{2}$ is the second Betti number of $\ov X$, and 
$\rho={\rm rk}\,\NS(\ov X)$. 
Since the ground field $k$ has characteristic zero, 
for almost all $\ell$ the group
$\H^3_{\et}(\ov X,\Z_\ell(1))$ is torsion-free: 
this is a consequence of the comparison theorem between \'etale cohomology
and Betti cohomology, see \cite[Thm. III.3.12]{milne}.
Thus the direct sum is finite.
By Hodge theory the condition $\H^2(X,O_{X})=0$ implies 
$\rho=b_{2}$. This shows that $\Br(\ov X)$ is finite.
As recalled above, the quotient $\Br_1(X)/\Br_0(X)$ is a subgroup
of $\H^1(K,\Pic(\ov X))$. 
The hypothesis $\H^1(X,O_{X})=0$ implies that the Picard variety of $X$
is trivial, hence $\Pic(\ov X) \simeq \NS(\ov X)$ is finitely generated.
In our case this group is torsion-free. It follows that
$\H^1(k,\Pic(\ov X))$ is finite, thus $\Br(X)/\Br_0(X)$ is also finite. 
This proves (i) and (ii).

When $\dim X=2$, for every prime $\ell$ 
there is a perfect duality pairing of finite abelian groups
$$\H^3_{\et}(\ov X,\Z_{\ell}(1))\{\ell\} \times \NS(\ov X)\{\ell\} \lra \Q_{\ell}/\Z_{\ell},$$
see \cite[III, 8.3]{G68}. Thus (iii) follows from (ii). \qed  

\section{Brauer group   of a variety over a local field and    Brauer pairing}

Throughout this section we use the following notation:
\medskip

$K$ is a finite extension of $\Q_p$,

$\ov K$ is an algebraic closure of $K$,

$K_{\rm nr}\subset \ov K$ is the maximal unramified extension,

$\g=\Gal(\ov K/K)$, \ $G=\Gal(K_{\rm nr}/K)$, \ $I=\Gal(\ov K/K_{\rm nr})$,

$R$ is the ring of integers of $K$,

$R_{\rm nr}$ is the ring of integers of $K_{\rm nr}$,

$X$ is a proper, smooth and geometrically integral variety over $K$,

$\pi:\X\to\Spec(R)$ is a faithfully flat proper morphism with
$X=\X\times_R K$,

$X_{\rm nr}=X\times_{K}K_{\rm nr}$,\ $\ov X=X\times_R \ov K$, \ 
$\X_{\rm nr}=\X\times_R {R_{\rm nr}}$. 

\medskip

\noindent We shall constantly use the classical result that
$\Br(K_{\rm nr})=0$, which implies
$\H^2(G,K_{\rm nr}^*)=\Br(K)$, see \cite{CL}, Ch. XII, Thm. 1 and its corollary.

The following finiteness statement  holds without a good reduction hypothesis.
It will not be used in this note.
 \bpr
For any smooth,  proper,   geometrically integral variety over $K$,
a finite field extension of $\Q_p$,
the group $\Ker[\Br(X) \to \Br(X_{\rm nr})]/\Br_0(X)$ is finite.
\epr
{\it Proof} There is a Hochschild--Serre spectral sequence attached 
to the morphism $X_{\rm nr}\to X$:
\begin{equation}
E_{2}^{pq}= \H^p(G,\H^q_{\et}(X_{\rm nr},\G_m))\Rightarrow \H^{p+q}_{\et}(X,\G_m).
\label{sss1}
\end{equation} 
Since $\H^2(G,K_{\rm nr}^*)=\Br(K)$, the exact sequence of low degree
terms of (\ref{sss1}) shows that the group under consideration
is a subgroup of $\H^1(G,\Pic(X_{\rm nr}))$. We have an exact sequence
of continuous discrete $\g$-modules
$$0\to \Pic^0(\ov X)\to\Pic(\ov X)\to \NS(\ov X)\to 0.$$
By the representability of the Picard functor over a field
of characteristic zero (due to Weil and Grothendieck,
see \cite{FAG}, Thm. 9.5.4 and Cor. 9.5.14) there exists an
abelian variety $A$ over $K$ such that $A(\ov K)$ is isomorphic
to $\Pic^0(\ov X)$ as a $\g$-module. Thus we rewrite the previous
sequence as
\begin{equation}
0\to A(\ov K)\to\Pic(\ov X)\to \NS(\ov X)\to 0. \label{ex1}
\end{equation}
The Hochschild--Serre spectral sequence attached 
to $\ov X\to X_{\rm nr}$ is 
$$E_{2}^{pq}=
\H^p(I,\H^q_{\et}(\ov X,\G_m))\Rightarrow \H^{p+q}_{\et}(X_{\rm nr},\G_m).
$$
By Hilbert's theorem 90 we have $\H^1(I,\ov K^*)=0$. Since $\Br(K_{\rm nr})=0$
we obtain that the natural map $\Pic(X_{\rm nr})\to \Pic(\ov X)^I$
is an isomorphism. Now, by taking $I$-invariants in (\ref{ex1}) 
we obtain the exact sequence of $G$-modules
$$0\to A(K_{\rm nr})\to\Pic(X_{\rm nr})\to \NS(\ov X)^I. $$
The group $\NS(\ov X)$ is finitely
generated by the theorem of N\'eron and Severi, hence so is
$\NS(\ov X)^I$. Thus there is a $G$-module $N$,
finitely generated as an abelian group, that fits into
the exact sequence of continuous discrete $G$-modules
$$0\to A(K_{\rm nr})\to\Pic(X_{\rm nr})\to N\to 0.$$
The resulting exact sequence of cohomology groups gives
an exact sequence
\begin{equation}
\H^1(G, A(K_{\rm nr}))\to\H^1(G,\Pic(X_{\rm nr}))\to 
\H^1(G,N). \label{ex2}
\end{equation}
We note that $G$ is canonically isomorphic to the profinite completion
$\hat \Z$, with the Frobenius as a topological generator.
If $M$ is a continuous discrete $G$-module which is finitely generated as
an abelian group, then $\H^1(G,M)$ is finite. To see this, let
$G'$ be a finite index subgroup of $G$ that acts trivially on $M$.
The group $G'\simeq \hat \Z$ has a dense subgroup $\Z$ generated by 
a power of the Frobenius.
Now $\H^1(G',M)$ is the group of continuous homomorphisms
$$\Hom_{\rm cont}(G',M)=\Hom_{\rm cont}(G',M_{\rm tors})=M_{\rm tors},$$
which is visibly finite. An application of the restriction-inflation sequence
finishes the proof of the finiteness of $\H^1(G,M)$.

To complete the proof of the proposition it remains to prove
the finiteness of $\H^1(G, A(K_{\rm nr}))$.
By Prop. I.3.8 of \cite{ADT} this group is isomorphic to
$\H^1(G, \pi_0(A_0))$,
where $\pi_0(A_0)$ is the group of connected components of the 
closed fibre $A_0$ of the N\'eron model of $A$ over $\Spec(R)$.
Since $\pi_0(A_0)$ is finite, we see that 
$\H^1(G, \pi_0(A_0))$ is finite. \qed

\medskip

Particular cases of the following two results
have been known for some time, see \cite{B07}.

\ble \label{pp0}
If $\X$ is smooth over $R$
with geometrically integral fibres, then

{\rm (i)} 
the following natural map is surjective:
$$\Br(K) \oplus \Ker[\Br(\X) \to \Br(\X_{\rm nr})]  \lra 
\Ker[\Br(X) \to \Br(X_{\rm nr})];$$

{\rm (ii)} for any $A \in \Ker[\Br(X) \to \Br(X_{\rm nr})] $
the image of the evaluation map $X(K) \to \Br(K)$ given by
$M \mapsto A(M)$ consists of one element.
\ele
{\it Proof} 
The map in (i) is well defined since $\Br(K_{\rm nr})=0$,
so that the composition $\Br(K)\to \Br(X)\to \Br(X_{\rm nr})$ is zero.

The restriction map $\Pic(\X_{\rm nr}) \to \Pic(X_{\rm nr})$ is surjective
since $\X_{\rm nr}$ is regular. 
The kernel of this map is generated by the classes of components 
of the closed fibre of $\X_{\rm nr}\to {\rm Spec}(R)$. 
The closed fibre is a principal divisor in $\X_{\rm nr}$;
since it is integral, the restriction map gives an isomorphism of $G$-modules
\begin{equation}
\Pic(\X_{\rm nr})\, \tilde\longrightarrow\, \Pic(X_{\rm nr}).
\label{iso2}
\end{equation}
There is a Hochschild--Serre spectral sequence attached 
to the morphism $\X_{\rm nr} \to \X$:
\begin{equation} E_{2}^{pq}=
\H^p(G,\H^q_{\et}(\X_{\rm nr},\G_m))\Rightarrow \H^{p+q}_{\et}(\X,\G_m),
\label{ss1}
\end{equation}
and a similar sequence (\ref{sss1}) attached to the morphism 
$X_{\rm nr}\to X$.
By functoriality the maps in (\ref{ss1}) and (\ref{sss1}) 
are compatible with the inclusion of the generic 
fibres $X\hookrightarrow \X$ and 
$X_{\rm nr}\hookrightarrow \X_{\rm nr}$. We have 
$\H^0_{\et}(\X_{\rm nr},\G_m)=R_{\rm nr}^*$ because $\pi:\X\to \Spec(R)$
is proper with geometrically integral fibres.
The low degree terms of the two spectral sequences
give rise to the following commutative diagram of exact sequences, 
where the equality is induced by (\ref{iso2}):
$$\begin{array}{ccccccc}
\H^2(G,R_{\rm nr}^*)&\to&\Ker[\Br(\X) \to \Br(\X_{\rm nr})]& 
\to &\H^1(G,\Pic(\X_{\rm nr})) &\to& \H^3(G,R_{\rm nr}^*)\\
\downarrow&&\downarrow&&||&&\downarrow\\
\H^2(G,K_{\rm nr}^*)&\to&\Ker[\Br(X) \to \Br(X_{\rm nr})]& \to &
\H^1(G,\Pic(X_{\rm nr})) &\to& \H^3(G,K_{\rm nr}^*)\\
\end{array}$$
We have $\H^3(G,R_{\rm nr}^*)=0$, since the group
$G\cong\hat{\mathbb Z}$ has
strict cohomological dimension 2. Since
$\H^2(G,K_{\rm nr}^*)=\Br(K)$, the
statement of (i) follows from the above diagram.

Any element $A \in \Br(\X) \subset \Br(X)$ vanishes
on $\X(R)=X(K)$, since it takes values in $\Br(R)=0$. This proves (ii).
\qed

\bpr \label{pp}
Assume that $\X$ is smooth over $R$
with geometrically integral fibres, $\H^1(X,O_{X})=0$ and 
the N\'eron--Severi group $\NS(\ov X)$ is torsion-free. 

{\rm (i)} 
Then the quotient $\Br_1(X)/\Br_0(X)$ is finite,
and every element of $\Br_1(X)\subset\Br(X)$ 
can be written as 
$\alpha+\beta$, where $\alpha\in\Br_0(X)$ and $\beta\in
\Br(\X)\subset\Br(X)$.

{\rm (ii)} Assume, moreover, that $\H^2(X,O_{X})=0$, and that
for every prime $\ell$
the group $\H^3_{\et}(\ov X,\Z_{\ell})$ is torsion-free.
Then the quotient $\Br(X)/\Br_0(X)$ is finite
and generated by the image of $\Br(\X)$.
\epr
{\it Proof} 
For $\ell\neq p$ the smooth base change theorem in 
\'etale cohomology for the smooth and proper morphism $\pi:\X \to \Spec(R)$
implies that the natural action of the inertia subgroup 
$I$ on $\H^2_{\et}(\ov X, \Z_{\ell}(1))$ is trivial.
Indeed, by \cite[Cor. VI.4.2]{milne} the \'etale sheaf $R^2\pi_*\mu_{\ell^m}$
is locally constant for every $m\geq 1$. Also, the fibre of $R^2\pi_*\mu_{\ell^m}$
at the generic geometric point $\Spec(\ov K)\to \Spec(R)$
is $\H^2_{\et}(\ov X, \mu_{\ell^m})$. 
Now it follows from Remark 1.2 (b) in \cite[Ch. V]{milne}
that the action of $\g$ on $\H^2_{\et}(\ov X, \mu_{\ell^m})$
factors through 
$$\pi_1(\Spec(R),\Spec(\ov K))=\Gal(K_{\rm nr}/K)=G=\g/I,$$
see \cite[Ex. I.5.2(b)]{milne}.
Thus $I$ acts trivially
on $\H^2(\ov X, \mu_{\ell^m})$ for every $m$, hence $I$ acts trivially
on $\H^2_{\et}(\ov X, \Z_{\ell}(1)).$

Since $K$ has characteristic zero, for any prime $\ell$ the Kummer sequence
gives a Galois equivariant embedding
$$\NS(\ov X) \otimes \Z_{\ell} \hookrightarrow \H^2_{\et}(\ov X, \Z_{\ell}(1)).$$
For any  $\ell\neq p$ we conclude that $I$ acts trivially on 
$\NS(\ov X) \otimes \Z_{\ell}$. 
Hence $I$ acts trivially on the quotient of the finitely generated 
group $\NS(\ov X)$ by its $p$-torsion subgroup.

So far the arguments apply to any smooth proper $R$-scheme
with geometrically integral fibres.

\medskip

If  $\NS(\ov X)$ is torsion-free, then we have 
$\NS(\ov X)\subset \NS(\ov X) \otimes \Z_{\ell}$, thus the natural 
action of $I$ on $\NS(\ov X)$ is also trivial.
Under the hypothesis $\H^1(X,O_{X})=0$, we now conclude that 
$\Pic(\ov X) \simeq \NS(\ov X)$ is a finitely
generated torsion-free abelian group with trivial action of inertia.
This implies $$\H^1(I,\Pic(\ov X))=0,$$
and $\H^1(k,\Pic(\ov X))$ is finite, hence $\Br_1(X)/\Br_0(X)$ is also finite.

We have the following commutative diagram of exact sequences
$$\begin{array}{ccccccc}
&&   && 0 &&0 \\
&&   &&\downarrow && \downarrow \\
0 & \to & \Ker[\Br(X) \to \Br(X_{\rm nr})] & \to &\H^2(G,K_{\rm nr}(X)^*) & \to & \H^2(G, \Div(X_{\rm nr}))\\
&& \downarrow &&\downarrow && \downarrow \\
0 & \to & \Ker[\Br(X) \to \Br(\ov X)]  & \to &\H^2(\g,\ov K(X)^*) & \to & \H^2(\g, \Div(\ov X))  \\
&&   &&\downarrow && \downarrow \\
&&   &&\H^2(I, \ov K(X)^*) & \to &\H^2(I, \Div(\ov X))     \\
\end{array}$$
The top and middle horizontal sequences are special cases of 
an exact sequence associated to a smooth variety over a field $K$ 
and a Galois extension of $K$, see \cite[\S 1.5.0]{CS}.
The two vertical sequences are restriction-inflation sequences. 
Their exactness
follows from $\H^1(I,\ov K(X)^*)=0$ (Hilbert's theorem 90) and
from $\H^1(I, \Div(\ov X))=0$ (since $X$ is smooth, the $I$-module 
$\Div(\ov X)$
is a permutation module, the vanishing then follows from Shapiro's lemma).

As mentioned above, under our hypotheses, $H^1(I,\Pic(\ov X))=0.$ 
We also have $\H^2(I,\ov K^*)=\Br(K_{\rm nr})=0$.
The exact sequence
$$0 \to \ov K^* \to \ov K(X)^* \to \Div(\ov X) \to  \Pic (\ov X) \to 0$$
then shows that the natural map 
$\H^2(I, \ov K(X)^*) \to \H^2(I, \Div(\ov X))$
is injective. From the commutative diagram above we conclude that the 
following natural inclusion is an isomorphism:
$$ \Ker[\Br(X)\to\Br(X_{\rm nr})] \tilde\lra\Ker[\Br(X)\to \Br(\ov X)].$$
An application of Lemma \ref{pp0} concludes the proof of statement (i).
 
Statement (ii) follows from (i) and Lemma  \ref{new1} (ii). \qed

\bigskip

We now would like to explore the situation when
$\H^2(X,O_{X})$ is not necessarily zero,
so we must  take  into account the transcendental Brauer group as well.
 
\bpr \label{p2}
Let $\ell$ be a prime, $\ell\not=p$.
Assume that $\X$ is smooth over $R$ with geometrically integral fibres, 
and that the closed geometric fibre  
has no connected unramified cyclic covering of degree $\ell$.
 
{\rm (i)} Then the group $\Br(X)\{\ell\}$
is generated by the images of
$\Br(\X)\{\ell\}$ and $\Br(K)\{\ell\}$.

{\rm (ii)} If $X(K)\neq \emptyset$, 
then for any  $A \in \Br(X)\{\ell\}$ the  image of
the evaluation map $X(K) \to \Br(K)$ given by
$M \mapsto A(M)$ consists of one element.
\epr
\emph{Proof} (i) Let $\F$ be the residue field of $R$,
and let $\X_0=\X\times_R \F$ be the closed fibre of $\pi$. 
By a special case of a result of K. Kato \cite[Prop. 1.7]{Kato}, 
to whose paper we refer for 
the explicit description of the maps involved, 
for any positive integer $n$ there is a natural {\it complex}
\begin{equation}
\Br(X)[{\ell^n}]
\buildrel{\res}\over{\hbox to 10mm{\rightarrowfill}}
\H^1(\F(\X_0),\Z/\ell^n) \lra \bigoplus_{Y\subset \X_0}\, \H^0(\F(Y),\Z/\ell^n(-1)),
\label{Kato}\end{equation}
where $Y$ ranges over closed  integral subvarieties of 
codimension 1 in $\X_0$, the field
$\F(\X_0)$ is the function field of $\X_0$,
and $\F(Y)$ is the function field of $Y$. 
The elements of $\H^1(\F(\X_0),\Z/\ell^n)$
correspond to characters of the absolute Galois group of $\F(\X_0)$
with values in $\Z/\ell^n$. These correspond to
connected cyclic coverings $W \to \X_0$ 
(not necessary unramified) of degree dividing $\ell^n$,
where $W$ is an irreducible normal variety. From Kato's complex it
follows that for any $A\in \Br(X)[{\ell^n}]$ the residue 
$\res(A)\in\H^1(\F(\X_0),\Z/\ell^n)$ 
is unramified in codimension 1.
Thus the corresponding covering $W \to \X_0$ is a covering
of a smooth scheme unramified at all the points of codimension 1,
hence is an \'etale covering of $\X_0$ by the 
Zariski--Nagata purity theorem (SGA 1 X, Th\'eor\`eme 3.1).
In other words, 
$$\res(A) \in \H^1_{\et}(\X_0,\Z/\ell^n) \subset     
\H^1(\F(\X_0),\Z/\ell^n).$$
Let $\ov\X_0=\X_0\times_{\F}\ov {\F}$, 
where $\ov {\F}$ is an algebraic closure of ${\F}$.
The spectral sequence
$$E_{2}^{pq}=\H^p(\F,\H^q_{\et}(\ov\X_0,\Z/\ell^n))\Rightarrow
\H^{p+q}_{\et}(\X_0,\Z/\ell^n)$$
gives rise to the exact sequence
$$0\to \H^1(\F,\Z/\ell^n)\to \H^1_{\et}(\X_0,\Z/\ell^n)\to
\H^1_{\et}(\ov\X_0,\Z/\ell^n).$$
Since $\ov\X_0$ has no connected unramified cyclic covering 
of degree $\ell$, we have $\H^1_{\et}(\ov\X_0,\Z/\ell^n)=0$.
Thus $\res(A)$ belongs to the injective image of $\H^1(\F,\Z/\ell^n)$
in $\H^1(\X_0,\Z/\ell^n)$.

By local class field theory the residue map 
$\Br(K)\{\ell^n\} \to \H^1(\F,\Z/\ell^n)$ is an isomorphism. Hence
for any $A \in \Br(X)\{\ell\}$ there exists $\alpha \in \Br(K)\{\ell^n\}$
such that the residue of $A-\alpha$ at any point of codimension 1 of 
$\X$ is zero. By Gabber's absolute purity theorem \cite{F},
 this implies that  $A-\alpha$ belongs to $\Br(\X)\{\ell\}
\subset \Br(X)\{\ell\}$. This completes the proof of (i).

Since $\X$ is proper over $\Spec(R)$ we have $X(K)=\X(R)$.
Statement (ii) follows because $\Br(R)=0$. \qed
 
\bigskip
 
\noindent{\bf Remarks} 
1. Already for $\pi$ smooth and proper,
it is an interesting $p$-adic problem to decide whether there 
is an analogous proposition for $\Br(X)\{p\}$.
For algebras split by an unramified extension of $K$, including
those of order divisible by $p$, this follows from
Lemma \ref{pp0} (see also  \cite[Prop. 6]{B07}). 

2. When $\X$ has
dimension at most 3, we may refer to Gabber's earlier purity theorem 
\cite[Thm. $2'$]{Ga} rather than to \cite{F}.

3. The hypotheses of Proposition \ref{p2} 
apply in particular  when the fibres 
of $\pi:\X \to \Spec(R)$ are smooth complete intersections of 
dimension at least 2 in projective space (an application of 
the weak Lefschetz theorem in \'etale cohomology, see \cite{Katz}).
In particular they apply to smooth surfaces of arbitrary degree in 
${\P}^3$.

\bigskip

\noindent{\bf A remark on the bad reduction case}.
Let $\ell$ be a prime, $\ell\not=p$.
Assume that $\X$ is a regular scheme, and $X(K)=\X(R)\neq \emptyset$. 
Let $\sZ$ be the largest open subscheme of $\X$ smooth over $\Spec(R)$,
such that
$\sZ\times_R K=X$, and every irreducible component of the closed fibre $\sZ_0$
is geometrically irreducible. Since $\X$ is regular,
a well known valuation argument shows that $X(K)=\X(R)=\sZ(R)$,
see, e.g., \cite{S96}, the proof of Lemma 1.1 (b).
Let $V_1,\ldots,V_n$ be the irreducible components
of $\sZ_0$, and let $\ov V_i=V_i\times_\F\ov \F$. Then
Kato's complex for $\sZ$ has the form
$$ \Br(X)[{\ell^n}]
\buildrel{\res}\over{\hbox to 10mm{\rightarrowfill}}
\bigoplus_{i=1}^n\H^1(\F(V_i),\Z/\ell^n) \lra 
\bigoplus_{Y\subset \sZ_0}\, \H^0(\F(Y),\Z/\ell^n(-1)).$$
We see that for any $A \in \Br(X)[\ell^n]$ the residue
$\res_i(A)\in \H^1(\F(V_i),\Z/\ell^n)$ belongs to the
subgroup $\H^1(V_i,\Z/\ell^n)$. We note that this group
is finite. This follows from the exact sequence
$$0\to\H^1(\F,\H^0(\ov V_i,\Z/\ell^n))\to
\H^1(V_i,\Z/\ell^n)\to \H^1(\ov V_i,\Z/\ell^n),$$
and the fact that $\H^1(\hat \Z,\Z/\ell^n)\simeq \Z/\ell^n$. 
We conclude that $\res_i$, as a function on $\Br(X)[\ell^n]$, takes only
finitely many values.

If the reduction $\tilde M$ of $M\in X(R)$ belongs to $V_i$,
then the local invariant of $A(M)\in \Br(K)$ is the pullback of
$\res_i(A)$ under the natural map $\H^1(V_i,\Z/\ell^n) \to
\H^1(\F,\Z/\ell^n)$ defined by $\tilde M$. Combining this information
for all irreducible components of $\sZ_0$ we obtain a finite partition
of the set $X(K)$ such that for any $A\in\Br(X)[\ell^n]$ the image of
the evaluation map $X(K) \to \Br(K)$ given by
$M \mapsto A(M)$ is constant on each element of the partition.

\section{Structure of the Brauer--Manin set over a number field}

We are now ready to prove the main results of this note.

\bthe \label{gen}
Let $k$ be a number field. Let $S$ be a finite set of places of $k$,
and let $\sO_{S}$ be the subring of $k$ consisting of the elements
that are integral at the primes not in $S$.
Let $\pi : \X \to \Spec(\sO_{S})$ be a smooth proper $\sO_S$-scheme 
with geometrically integral fibres.
Let $X/k$ be its generic fibre. Assume 

{\rm (i)}  $\H^1(X,O_{X})=0$;

{\rm (ii)} the N\'eron--Severi group $\NS(\ov X)$ has no torsion;

{\rm (iii)} $\Br(X)/\Br_1(X)$ is a finite abelian group
of order invertible in $\sO_{S}$.

\noindent Then the answer to Question $1$ is in the affirmative.
\ethe
{\it Proof} Let $v$ be a place of $k$ not contained in $S$, and let
$p$ be the residual characteristic of $k_v$. Assumption (iii)
implies that $\Br(X)\{p\}\subset \Br_1(X)$. Then, by Proposition \ref{pp},
the image of $\Br(X)\{p\}$ in $\Br(X_v)$ is contained in the subgroup
generated by the images of $\Br(k_v)$ and $\Br(\X_v)$.
By Proposition \ref{p2} the same is true for $\Br(X)\{\ell\}$ for 
any prime $\ell$ not equal to $p$. For this we only need to check
that $\H^1_{\et}(\ov\X_M,\Z/\ell)=0$, where $\ov\X_M$ 
is the closed geometric fibre
of $\pi:\X_v\to \Spec(\sO_v)$. By the smooth base change theorem
for \'etale cohomology (see, e.g. \cite{milne}, VI, Cor. 4.2)
the group $\H^1_{\et}(\ov\X_M,\Z/\ell)$ is isomorphic to 
$\H^1_{\et}(\ov X_v,\Z/\ell)$, which in turn is isomorphic to 
$\H^1_{\et}(\ov X,\Z/\ell)$ by \cite{milne}, VI, Cor. 4.3.
The Kummer exact sequence gives an isomorphism
$\H^1_{\et}(\ov X,\mu_\ell)\tilde\lra\Pic(\ov X)[\ell]$,
and the vanishing of the latter group follows 
from conditions (i) and (ii). 
The statement now follows from Lemma \ref{l2}. \qed

\bco\label{t1}
Let $\pi : \X \to \Spec(\sO_{S})$ be a smooth proper $\sO_S$-scheme 
with geometrically integral fibres.
Let $X/k$ be its generic fibre. Assume 

{\rm (i)}  $\H^{i}(X,O_{X})=0$ for $i=1,\,2$;

{\rm (ii)} the N\'eron--Severi group $\NS(\ov X)$ has no torsion;

{\rm (iii)} either $\dim X=2$, or
$\H^3_{\et}(\ov X,\Z_{\ell})$ is torsion-free for every prime $\ell$
outside $S$.

\noindent Then the answer to Question $1$ is in the affirmative.
\eco
{\it Proof} This follows from Theorem \ref{gen} by
Lemma \ref{new1} and its proof. \qed

\medskip

This corollary can be applied to rationally connected varieties.
Indeed, over a
 field of characteristic zero, these varieties
are
$O_X$-acyclic and algebraically simply connected
\cite[Cor. 4.18]{debarre}).

When we no longer have $\H^2(X,O_{X})=0$,  
condition (iii) in Theorem \ref{gen}  is not easy to check in general.
However, this can be done in an important particular case:
as an application of rather delicate computations of \cite{ISZ}
we now show that the answer to Question 1 is positive for
smooth diagonal quartics over $\Q$, so that only the real place and the primes
of bad reduction can show up in the Brauer--Manin obstruction.

\bco \label{t2}
Let $D$ be the diagonal quartic surface over $\Q$ given by
\begin{equation}
x_0^4+a_1x_1^4+a_2x_2^4+a_3x_3^4=0, \label{1}
\end{equation}
where $a_1,\,a_2,\,a_3\in \Q^*$.
Let $\sS$ be the set of primes consisting of $2$ and the primes
dividing the numerators or the denominators
of $a_1,\,a_2,\,a_3$.
Let $Z$ be the image of the projection
$$D(\A_\Q)^\Br\rightarrow D(\R)\times\prod_{p\in \sS}D(\Q_p).$$
Then $D(\A_\Q)^\Br=Z\times \prod_{p\not\in \sS}D(\Q_p)$.
\eco
\emph{Proof} $D$ is a K3 surface, for which the geometric conditions
(i) and (ii) of Theorem \ref{gen} are well known.
By Theorem 3.2 of \cite{ISZ}, 
only the primes from $\{2,\,3,\,5\}\cap \sS$
can divide the order of the finite group $\Br(D)/\Br_1(D)$. \qed

\medskip

There are other K3 surfaces over $\Q$ to which Theorem \ref{gen} can be applied.
Let $X$ be the Kummer surface attached to the product of elliptic curves
$E$ and $E'$ over $\Q$. By \cite[Prop. 1.4]{SZ2} we have 
$$\Br(\ov X)^\Gamma\cong \Br(\ov E\times \ov E')^\Gamma,$$
where $\Gamma=\Gal(\ov\Q/\Q)$. If
$E$ and $E'$ are not isogenous over $\ov \Q$, then by \cite[Prop. 3.1]{SZ2}
for any integer $n$ we have
$$\Br(\ov E\times \ov E')[n]^\Gamma\cong \Hom_\Gamma(E[n],E'[n]).$$
In \cite{SZ2} (Prop. 4.2 and 
Example A3) we construct infinitely many pairs of non-isogenous 
elliptic curves
$E$, $E'$ such that $\Hom_\Gamma(E[\ell],E'[\ell])=0$ for any odd prime
$\ell$. Then $\Br(X)/\Br_1(X)$ is a finite abelian 2-group, so that
only the archimedean place, $2$ and the primes of bad reduction
of $E$ and $E'$ can turn up in the description of the Brauer--Manin
set of $X$.

\section{Transcendental Brauer group and the Tate conjecture for divisors} \label{Tate}

In connection with Lemma \ref{l1} let us briefly discuss
the following  question (see also \cite{SZ}).

\medskip

\noindent{\bf Question 2}. Let $X$ be a smooth,
proper and geometrically integral variety 
over a field $k$ finitely generated over
$\Q$. We have the inclusion of groups  
$\Br(X)/\Br_1(X) \subset \Br(\ov X)^\Gamma$. Are these
two groups finite?

\medskip

It is well known that the finiteness of $\Br(\ov X)^\Gamma$
is related to the Tate conjecture for divisors. This conjecture says
that if $k$ is a field finitely generated over $\Q$, then
for any prime $\ell$ the natural inclusion
$$(\NS(\ov X)\otimes_{\Z}\Q_\ell)^\Gamma\hookrightarrow
\H^2_{\et}(\ov X,\Q_\ell(1))^\Gamma$$
should be an isomorphism.
The following partial answer to Question 2 must have been known 
to many people. It was noticed by one of us ten years ago.

\bpr\label{f}
Let $X$ be a smooth, proper and geometrically integral variety 
over a field $k$ finitely generated over
$\Q$. Assume the $\ell$-adic Tate conjecture for divisors.

{\rm (i)} If $\dim X=2$, then $\Br(\ov X)^\Gamma \{\ell\}$ is finite.

{\rm (ii)} For $X$ of any dimension assume in addition the
semisimplicity of the continuous
$\Gamma$-module $\H^2_{\et}(\ov X,\Q_\ell(1))$. 

Then
$\Br(\ov X)^\Gamma \{\ell\}$ is finite.
\epr
{\it Proof} 
   Let $T_\ell(\Br(\ov X))$ be the $\ell$-adic Tate module
of $\Br(\ov X)$, and let
$V_\ell(\Br(\ov X))=T_\ell(\Br(\ov X))\otimes_{\Z_\ell}\Q_\ell$.
The Kummer sequence gives a well known exact sequence
of continuous $\Gamma$-modules
$$
0\lra \NS(\ov X)\otimes_\Z \Q_\ell \lra \H^2_{\et}(\ov X,\Q_\ell(1))
\lra V_\ell(\Br(\ov X)) \lra 0.
$$
When $X$ is a surface, the cup-product defines a 
non-degenerate Galois-equivariant bilinear pairing on $\H^2_{\et}(\ov X,\Q_\ell(1))$
with values in $\Q_\ell$.
The restriction of this pairing to $\NS(\ov X)\otimes_\Z \Q_\ell$
is non-degenerate, hence we obtain a direct sum decomposition of 
$\Gamma$-modules
$$\H^2_{\et}(\ov X,\Q_\ell(1))\cong 
(\NS(\ov X)\otimes_\Z \Q_\ell)\oplus V_\ell(\Br(\ov X)).$$
We have the same conclusion under the semisimplicity  
hypothesis of part (ii).
Now Tate's conjecture implies $V_\ell(\Br(\ov X))^\Gamma=0$.

For any abelian group $A$, the Tate module $T_{\ell}(A)$ is
a torsion-free $\Z_{\ell}$-module. If $A$ is a torsion group
whose $\ell$-primary component is of cofinite type, 
  $T_{\ell}(A)\otimes_{\Z_\ell}\Q_\ell/\Z_\ell$
is the maximal $\ell$-divisible subgroup of $A$.
 We thus have an exact sequence
\begin{equation}
0\lra T_\ell(\Br(\ov X))\lra V_\ell(\Br(\ov X))\lra
T_\ell(\Br(\ov X))\otimes_{\Z_\ell}\Q_\ell/\Z_\ell \lra 0.
\label{seq1}
\end{equation}
where $T_\ell(\Br(\ov X))\otimes_{\Z_\ell}\Q_\ell/\Z_\ell$
is the maximal divisible subgroup $\Br(\ov X)\{\ell\}_{\rm div}$ of 
$\Br(\ov X)\{\ell\}$.
The exact sequence of Galois cohomology attached to
(\ref{seq1}) gives an exact sequence:
\begin{equation}
V_\ell(\Br(\ov X))^\Gamma\lra (\Br(\ov X)\{\ell\}_{\rm div})^\Gamma
\lra \H^1(\Gamma, T_\ell(\Br(\ov X))).
\label{seq2}
\end{equation}
By a general result of  \cite[Prop. 2.3]{Tate} the kernel of the last arrow in (\ref{seq2})
is the maximal divisible subgroup of $(\Br(\ov X)\{\ell\}_{\rm div})^\Gamma$.
From $V_\ell(\Br(\ov X))^\Gamma=0$ we conclude that
 the maximal divisible subgroup of $(\Br(\ov X)\{\ell\}_{\rm div})^\Gamma$
is zero.
Since $(\Br(\ov X)\{\ell\}_{\rm div})^\Gamma$ is a subgroup of
$\Br(\ov X)\{\ell\}$, which is a torsion group of finite cotype,
we see that $(\Br(\ov X)\{\ell\}_{\rm div})^\Gamma$ is also 
a torsion group of finite cotype. It follows that
$(\Br(\ov X)\{\ell\}_{\rm div})^\Gamma$ is finite.

As already recalled, by \cite{G68}, III (8.9), 
we have an exact sequence of continuous
$\Gamma$-modules
$$0\lra \Br(\ov X)\{\ell\}_{\rm div} \lra \Br(\ov X)\{\ell\}
\lra \H^3_{\et}(\ov X,\Z_\ell(1))_{\rm tors}\lra 0.$$
The attached exact sequence of Galois cohomology gives
an exact sequence
$$0\lra (\Br(\ov X)\{\ell\}_{\rm div})^\Gamma \lra \Br(\ov X)\{\ell\}^\Gamma
\lra \H^3_{\et}(\ov X,\Z_\ell(1))_{\rm tors}^\Gamma$$
Since $\H^3_{\et}(\ov X,\Z_\ell(1))_{\rm tors}$ is finite, our
statement follows. \qed

\bigskip

\noindent{\bf Remarks} 1. When $\H^2(X,O_{X})=0$ the group
$\Br(\ov X)$ is finite, so Question 2 is trivial.

2. For abelian varieties or K3 surfaces the Tate conjecture 
for divisors is known. Here one can do better that Proposition
\ref{f}: by the main theorem of \cite{SZ}
Question 2 has a positive answer in these cases.

\bigskip

\noindent
CNRS, UMR 8628, Math\'ematiques, B\^atiment 425, Universit\'e Paris-Sud,
F-91405 Orsay, France
\medskip

\noindent jlct@math.u-psud.fr

\bigskip

\bigskip

\noindent Department of Mathematics, South Kensington Campus,
Imperial College London,
SW7 2BZ England, U.K.

\smallskip

\noindent Institute for the Information Transmission Problems,
Russian Academy of Sciences, 19 Bolshoi Karetnyi,
Moscow, 127994 Russia
\medskip

\noindent a.skorobogatov@imperial.ac.uk

\end{document}